\documentclass{amsart}

\usepackage{amsmath}
\usepackage{amsfonts}
\usepackage{amssymb}

\numberwithin{equation}{section}

\newtheorem{theorem}[subsection]{Theorem}

\newtheorem{lemma}[subsection]{Lemma}

\newcommand{\la}{\lambda}
\newcommand{\ep}{\epsilon}
\newcommand{\ii}{{\bf i}}
\newcommand{\sumn}{\frac{1}{x}\sum_{n\le x}}
\newcommand{\prodi}{\prod_{i=1}^{r}}
\newcommand{\prodp}{\prod_{y<p<z}}
\newcommand{\sump}{\sum_{y<p<z}}
\newcommand{\om}{\omega}
\newcommand{\sumi}{\sum_{i=1}^{r}}

\newcommand{\al}{\alpha}

\newcommand{\de}{\delta}

\newcommand{\half}{\frac{1}{2}}

\begin{document}

\title{Spacings between integers having typically many prime factors}
\author{Rizwanur Khan}
\address{University of Michigan, Department of Mathematics,
         530 Church St., Ann Arbor, MI 48109, USA}
\email{rrkhan@umich.edu}

\begin{abstract}
We show that the sequence of integers which have nearly the typical
number of distinct prime factors forms a Poisson process. More
precisely, for $\de$ arbitrarily small and positive, the nearest
neighbor spacings between integers $n$ with $|\om(n)-\log_2 n|\le
(\log_2 n)^{\de}$ obey the Poisson distribution law.
\end{abstract}

\maketitle

\section{Introduction}

Consider $n$ random variables independently and uniformly taking
real values in the interval $[0,n]$. Let $Y_1<...<Y_n$ denote the
order statistics obtained by arranging these random variables in
increasing order. Setting $Y_0=0$ and $Y_{n+1}=n$, let
$D_i=Y_{i+1}-Y_{i}$ for $0\le i \le n$ denote the nearest neighbor
spacings of the order statistics. Thus $D_1+\ldots +D_n=n$ and by
symmetry it follows that for $0<\la<n$ a real number,
\begin{align*}
\text{Prob}(D_i>\la)=\text{Prob}(D_1>\la)=\Big(\frac{n-\la}{n}\Big)^n.
\end{align*}
Thus $\text{Prob}(D_i>\la)\sim e^{-\la}$ as $n\rightarrow \infty$.
This is the exponential or Poisson distribution.

We are interested in the spacing distributions of arithmetic
sequences. An example of such a sequence is the sequence of prime
numbers less than $x$, which form a sparse subset of the integers of
density $1/\log x$ by the prime number theorem. This is similar to
$Y_1,\ldots,Y_n$ being sparse in the interval $[0,n]$. If $p_i$
denotes the $i$-th prime less than $x$, we rescale to consider
instead the sequence $\widetilde{p}_i=p_i/\log x$ of `normalized'
primes so that that the average spacing between consecutive
normalized primes is $1$ as $x\rightarrow \infty$. This matches the
expected value of $D_i$ above as $n\rightarrow \infty$. Gallagher
\cite{gall} showed that assuming the validity of the
Hardy-Littlewood prime $k$-tuple conjectures, we have for $\la
>0$ real that
\begin{align}
\frac{1}{x} \# \{ i\le x:  \widetilde{p}_{i+1}-\widetilde{p}_i
> \lambda \}\sim
 e^{-\la}
\end{align}
as $x\rightarrow \infty$. Thus conditionally we see that the
spacings between primes obey the Poisson distribution law, as in the
prototypical situation of randomly dispersed objects mentioned at
the start. More recently Kurlberg and Rudnick \cite{kurlrudn} showed
that the spacings between quadratic residues modulo $q$ , as the
number of distinct prime divisors of $q$ tends to infinity, follow
the Poisson distribution. There are many other interesting
arithmetic sequences that are conjectured to be Poisson processes,
but only few examples exist with proof. For example, it is an open
problem to show that the spacings between the fractional parts of
$n^2 \sqrt{2}$ for $n\le x$, as $x\rightarrow \infty$ are Poisson
distributed (see \cite{sarnak}). The reader may find a few more
examples of such work listed in the references section (see
\cite{cob, grankurl, hool}). Of course there are important
arithmetic sequences which are not expected to behave like randomly
dispersed elements in this sense, such as the non-trivial zeros of
the Riemann Zeta function. In this paper we are interested in the
spacings between integers with not only one prime factor as in
Gallagher's work, but with the typical number of distinct prime
factors. We first explain what is meant by `typical'.

Let $\omega(n)$ denote the number of distinct prime factors of $n$.
It is easy to see that integers $n\le x$ have $\log\log x$ distinct
prime factors on average:
\begin{align}
\label{first}\sumn \om(n) = \sumn \sum_{p|n} 1=
\frac{1}{x}\sum_{p\le x} \sum_{\substack{n\le x\\ p|n}} 1 =
\frac{1}{x} \sum_{p\le x} \Big\lfloor \frac{x}{p} \Big\rfloor =
\log_2 x +O(1),
\end{align}
where we write $\log_2 x$ for $\log\log x$, and similarly for
$\log_j x$. Also, throughout this paper $p$ and $q$ will be used to
denote primes. The variance can be shown to be
\begin{align}
\label{variance}\sumn (\om(n)-\log_2 x)^2\sim \log_2 x.
\end{align}
Note that (\ref{first}) and (\ref{variance}) imply that $\om(n)\sim
\log_2 x$ for all most all $n\le x$. Erd{\H o}s and Kac \cite{erdo}
further showed that $\om(n)$ is normally distributed with mean
$\log_2 x$ and standard deviation $\sqrt{\log_2 x}$. R\'{e}nyi and
Tur\'{a}n \cite{reny} proved this result with a sharp error term.
The following theorem can also be found in Tenenbaum's book
\cite{tenen}.
\begin{theorem} \label{normal} Given a real number $C>0$ we have for $0<c<C$ that the number of integers $n\le x$ for which
$-c<\frac{\omega(n)-\log_2x}{\sqrt{\log_2 x}}<c$ is
\begin{align*}
x\frac{1}{\sqrt{2\pi}}\int_{-c}^{c}\exp(-u^2/2)du
+O_C(x/\sqrt{\log_2 x}).
\end{align*}
\end{theorem}
\noindent We \cite{khan} proved a slightly weaker version of the
Theorem \ref{normal} by methods similar to those in this paper.

We conjecture that the spacings between integers $n\le x$ with
$|\om(n)-\log_2 x|\le \sqrt{\frac{\pi}{2}}$ (that is, integers with
more or less exactly $\log_2 x$ distinct prime factors) obey the
Poisson distribution law but we are unable to prove it. Instead we
look at an easier question. For any fixed $0<\de<1/2$, let us say an
integer less than $x^2$ is `$\de$-normal' if 
\begin{align*}
|\om(n)-\log_2 x|\le
\sqrt{\frac{\pi}{2}} (\log_2 x)^{\de}.
\end{align*}
We study the sequence of
$\de$-normal numbers. These are integers having nearly the expected
number of prime factors, as $(\log_2 x)^{\de}$ is smaller than the
standard deviation $\sqrt{\log_2 x}$ of $\om(n)$. Denote the
sequence of $\de$-normal numbers in increasing order by
$N_1,N_2,N_3\ldots$. Up to $x$, there are $x (\log_2 x)^{-1/2+\de}$
such integers by Theorem \ref{normal}, since an integer is
$\de$-normal if and only if $\Big|\frac{\om(n)-\log_2
x}{\sqrt{\log_2 x}}\Big|\le \sqrt{\frac{\pi}{2}}(\log_2
x)^{-1/2+\de}.$ Thus we should rescale these integers by setting
$\widetilde{N}_i={N}_i (\log_2 x)^{-1/2+\de}$. Our main theorem is
\begin{theorem}  \label{main} For any fixed real number $\lambda>0$ we
have
\begin{align*} \frac{1}{x} \# \{ i\le x:  \widetilde{N}_{i+1}-\widetilde{N}_i
> \lambda \}\sim
 e^{-\la}.
\end{align*}
\end{theorem}
\noindent Throughout this paper, all implicit constants may depend
on $\de$ and $\la$.

\section{Independence between additive shifts of the $\om(n)$ function}

In this section we show how Theorem \ref{main} can be reduced to
studying correlations between the additive shifts of the function
$\om(n)$. We will show for example that $\om(n)-\log_2 x$,
$\om(n+1)-\log_2 x$, and $\om(n+2)-\log_2 x$ behave independently.
Define ${\mathcal N}(x)$ to be the number of $\de$-normal integers less than $x$. The left hand side of Theorem \ref{main}
is asymptotic to
\begin{align}
\nonumber &\frac{1}{x} \# \{ i\le x: N_{i+1}-N_i
> \lambda (\log_2 x)^{1/2-\de}  \} \\
\label{interval} &\sim \frac{1}{x}\# \{ N\le x (\log_2 x)^{1/2-\de}
: {\mathcal N}\big(N+\lambda (\log_2 x)^{1/2-\de}\big)-{\mathcal
N}\big(N\big) =0 \},
\end{align}
where $N$ denotes a $\de$-normal number. Define ${\mathcal
N}_{b_1,...,b_r}(x)$ to be the number of integers $n\le x$ for which
$n+b_i$ is $\de$-normal for all $1\le i \le r$ and let $\sigma(m,r)$
denote the number of maps from the set $\{1,...,m\}$ onto
$\{1,...,r\}$. We have the $m$-th moment of ${\mathcal
N}\big(N+\lambda (\log_2 x)^{1/2-\de}\big)-{\mathcal N}\big(N\big)$:
\begin{align}
&\nonumber \frac{1}{x} \sum_{N\le x(\log_2 x)^{1/2-\de}}
\Big({\mathcal N}\big(N+\lambda
(\log_2 x)^{1/2-\de}\big)-{\mathcal N}\big(N\big)\Big)^{m}\\
&\label{moment}=\frac{1}{x}\sum_{r=1}^{m}\sigma(m,r) \sum_{1\le
b_1<...<b_r\le \lambda (\log_2 x)^{1/2-\de}} {\mathcal
N}_{0,b_1,...,b_r}\big(x(\log_2 x)^{1/2-\de}\big).
\end{align}
We will prove
\begin{theorem} \label{prop}
For a fixed integer $r$ and any integers $0\le b_1<...<b_r\le
\lambda (\log_2 x)^{1/2-\de}$, we have
\begin{align*}
\frac{1}{x}{\mathcal N}_{b_1,...,b_r}(x) \sim (\log_2
x)^{(-1/2+\de)r}.
\end{align*}
\end{theorem}
\noindent Throughout this paper all implicit constants may depend
$r$. Since a randomly chosen integer less than $x$ is $\de$-normal
with probability $(\log_2 x)^{-1/2+\de}$, the theorem above says
that $n+b_1,\ldots,n+b_r$ are independently likely to be
$\de$-normal. Theorem \ref{prop} implies that for fixed $m$ we have
that (\ref{moment}) is asymptotic to
\begin{align}
\label{pmoment}\sim \sum_{r=1}^{m}\sigma(m,r)\frac{\lambda^r}{r!}=
\sum_{j=0}^{\infty} j^m \frac{e^{-\la}\la^j}{j!},
\end{align}
which is the $m$-th moment of the Poisson distribution (the identity
above is known as Dobinski's formula). The Poisson distribution can
be recovered from these moments. Let us sketch this; we have that (\ref{interval}) is
\begin{align}
&\nonumber\frac{1}{x}\# \{ N\le x (\log_2 x)^{1/2-\de} : {\mathcal
N}\big(N+\lambda (\log_2 x)^{1/2-\de}\big)-{\mathcal N}\big(N\big)=0
\}\\
&\nonumber\sim1-\frac{1}{x}\sum_{j=1}^{\infty} \sum_{\substack{N\le
x (\log_2 x)^{1/2-\de}\\{\mathcal N}\big(N+\lambda (\log_2
x)^{1/2-\de}\big)-{\mathcal
N}\big(N\big)=j}} 1\\
&\label{recover}=1-\frac{1}{x}\sum_{j=1}^{\infty}
\sum_{m=0}^{\infty} \frac{(2\pi \ii)^m}{m!}j^m \sum_{\substack{N\le
x (\log_2 x)^{1/2-\de}\\{\mathcal N}\big(N+\lambda (\log_2
x)^{1/2-\de}\big)-{\mathcal N}\big(N\big)=j}} 1.
\end{align}
Now
\begin{align*}
\frac{1}{x}\sum_{j=1}^{\infty} j^m \sum_{\substack{N\le x (\log_2
x)^{1/2-\de}\\{\mathcal N}\big(N+\lambda (\log_2
x)^{1/2-\de}\big)-{\mathcal N}\big(N\big)=j}} 1
\end{align*}
is the $m$-th moment of ${\mathcal N}\big(N+\lambda (\log_2
x)^{1/2-\de}\big)-{\mathcal N}\big(N\big)$. By (\ref{pmoment}) we
get (an explicit dependence on $m$ of the error term is not needed)
that (\ref{recover}) is asymptotic to
\begin{align}
\sim 1-\sum_{j=1}^{\infty} \sum_{m=0}^{\infty} \frac{(2\pi
\ii)^m}{m!}j^m \frac{e^{-\la}\la^j}{j!}=e^{-\la}.
\end{align}
Thus Theorem \ref{main} follows from Theorem \ref{prop}. Next we discuss the demonstration
of Theorem \ref{prop}.

The characteristic function of a random variable with a normal
distribution is $\exp(-T^2/2)$. We show the independence of
$\frac{\om(n+b_i)-\log_2 x}{\sqrt{\log_2 x}}$ for $1\le i \le r$ by
showing that their joint characteristic function equals essentially
$\prodi \exp(-T_i^2/2)$. Actually it is more convenient to work with
$\om(n;y,z)$ in place of $\om(n)$, where we set
\begin{align*}
y=y(x)=(\log x)^{3r}
\end{align*}
and
\begin{align*}
z=z(x)=x^{((\log_2 x)^{-3r})}
\end{align*}
and define
\begin{align*}
\om(n;y,z)=\sum_{\substack{p|n\\y<p<z}}1.
\end{align*}
Accordingly we work with $\om(n;y,z)-\sump \frac{1}{p}$ in place of
$\om(n)-\log_2 x$. We will soon see that there is not much loss in
disregarding the primes less than $y$ or greater than $z$. In an
imprecise sense, the reason for this is that on average integers
have few small prime factors and few large prime factors. In the
next section we will prove the following theorem.
\begin{theorem} \label{fou}
Let $t_i=T_i(\sump \frac{1}{p})^{-1/2}$ be real. For $|T_i| \le
\frac{1}{1000}(\sump \frac{1}{p})^{1/2}$, we have
\begin{align*}
\sumn \prodi \exp \Big(\ii T_i \frac{\om(n+b_i;y,z)-\sump
\frac{1}{p}}{\sqrt{\sump \frac{1}{p}}}\Big)= \prodi &\exp\Big( \Big(
e^{\ii t_i}-1-\ii t_i \Big) \sump \frac{1}{p}\Big)\\
&+O(1/\log x).
\end{align*}
\end{theorem}
\noindent Observe that for $T_i\le (\log_2 x)^{\ep}$ for small
enough $\ep>0$ we have
\begin{align}
\nonumber \exp\Big( \Big( e^{\ii t_i}-1-\ii t_i \Big) \sump
\frac{1}{p}\Big)&=\exp\Big(\Big(\frac{-t_i^2}{2}+O(t_i^3)\Big)\sump
\frac{1}{p}\Big)\\
&\label{observe}=\exp\Big(\frac{-T_i^2}{2}\Big)\Big(1+O\Big(\frac{T_i^3}{\sqrt{\log_2
x}}\Big)\Big),
\end{align}
and for $(\log_2 x)^{\ep}<T_i\le (\log_2 x)^{1/2-\ep}$ we have
\begin{align}
\label{observe2} \exp\Big( \Big( e^{\ii t_i}-1-\ii t_i \Big) \sump
\frac{1}{p}\Big)\ll \exp(-T_i^2/4),
\end{align}
where the implied constants depend on $\ep$.

To see how Theorem \ref{fou} implies Theorem \ref{prop} we will use
the following lemmas.
\begin{lemma}\label{inversion}
Let $\psi(x)$ be a real function differentiable $\lfloor 4r/\de
\rfloor$ times and satisfying
\begin{align}
&\nonumber 0\le \psi(x) \le 1 \text{ for } x\in \mathbb{R},\\
&\nonumber \psi(x)=0 \text{ for } |x| \ge 2(\log_2
x)^{-1/2+\de}\\
&\nonumber \int_{-\infty}^{\infty} \psi(x) dx\sim \sqrt{2\pi}(\log_2
x)^{-1/2+\de},\\
&\label{psi}|\psi^{(j)}(x)|\ll (\log_2 x)^{j(1-\de)/2} \text { for
any positive integer $j\le \lfloor 4r/\de \rfloor$}.
\end{align}
We have
\begin{align*}
\sumn \prodi \psi\Big(\frac{\om(n+b_i;y,z)-\sump
\frac{1}{p}}{\sqrt{\sump \frac{1}{p}}}\Big)\sim (\log_2
x)^{(-1/2+\de)r}.
\end{align*}
\end{lemma}
\proof
Let
\begin{align*}
\hat{\psi}(T)=\int_{-\infty}^{\infty} \psi(u)e^{-\ii uT}du
\end{align*}
 denote the Fourier transform
of $\psi$. By Fourier inversion we have
\begin{align}
\nonumber &\sumn \prodi \psi\Big(\frac{\om(n+b_i;y,z)-\sump
\frac{1}{p}}{\sqrt{\sump \frac{1}{p}}}\Big)\\
&\label{inv1}=\sumn \prodi \frac{1}{2\pi} \int_{-\infty}^{\infty}
\hat{\psi}(T_i) \exp \Big(\ii T_i \frac{\om(n+b_i;y,z)-\sump
\frac{1}{p}}{\sqrt{\sump \frac{1}{p}}}\Big)dT_i.
\end{align}
We have that $|\hat{\psi}(T_i)|\ll (\log_2
x)^{j(1-\de)/2}|T_i|^{-j}$, by integrating by parts $j$ times and
using (\ref{psi}). Thus (\ref{inv1}) equals
\begin{align}
\sumn &\prodi \frac{1}{2\pi} \int_{-(\log_2 x)^{1/2-\de/4}}^{(\log_2
x)^{1/2-\de/4}} \hat{\psi}(T_i) \exp \Big(\ii T_i
\frac{\om(n+b_i;y,z)-\sump \frac{1}{p}}{\sqrt{\sump
\frac{1}{p}}}\Big)dT_i\\
&\nonumber+O\big((\log_2 x)^{-r}\big).
\end{align}
Now by Theorem \ref{fou} and observations (\ref{observe}) and
(\ref{observe2}), the main term above equals
\begin{align}
&\nonumber \prodi \Bigg(\frac{1}{2\pi} \int_{-(\log_2
x)^{\ep}}^{(\log_2 x)^{\ep}} \hat{\psi}(T_i)
\exp\Big(\frac{-T_i^2}{2}\Big)dT_i +O\Big(\frac{1}{\sqrt{\log_2
x}}\Big)\Bigg)\\
&\label{inv2}= \Bigg(\frac{1}{2\pi} \int_{-\infty}^{\infty}
\hat{\psi}(T) \exp\Big(\frac{-T^2}{2}\Big)dT
+O\Big(\frac{1}{\sqrt{\log_2 x}}\Big)\Bigg)^r.
\end{align}
Recall that the Fourier transform of $\tfrac{1}{\sqrt{2\pi}
}\exp\big(\frac{-u^2}{2}\big)$ is $\exp \big(\frac{-T^2}{2}\big)$.
By the Plancherel formula, (\ref{inv2}) equals
\begin{align}
\Bigg(\frac{1}{\sqrt{2\pi}} \int_{-\infty}^{\infty} \psi(u)\exp
\Big(\frac{-u^2}{2}\Big)du+O\Big(\frac{1}{\sqrt{\log_2
x}}\Big)\Bigg)^r \sim (\log_2 x)^{(-1/2+\de)r}.
\end{align}
\endproof
\noindent To prove Theorem \ref{prop} we need to show
\begin{align*}
\sumn \prodi \psi \Big(\frac{\om(n+b_i)-\log_2 x}{\sqrt{\log_2
x}}\Big)\sim (\log_2 x)^{(-1/2+\de)r},
\end{align*}
where $\psi$ is a suitable smooth function approximating the
characteristic function of the interval
$[-\sqrt{\frac{\pi}{2}}(\log_2
x)^{-1/2+\de},\sqrt{\frac{\pi}{2}}(\log_2 x)^{-1/2+\de}]$. This is
accomplished by Lemma \ref{inversion}, provided that we can show
that we may neglect prime factors smaller than $y$ or larger than
$z$ without significant loss. This is the purpose of the next lemma.
\begin{lemma}\label{notloss} Except for $O\big(x(\log_2 x)^{-r}\big)$ integers less
than $x$ we have
\begin{align*}
\Bigg|\frac{\om(n)-\log_2 x}{\sqrt{\log_2 x}}-\frac{\om(n;y,z)-\sump
\frac{1}{p}}{\sqrt{\sump \frac{1}{p}}}\Bigg|\ll (\log_2
x)^{-1/2+\de/2}.
\end{align*}
\end{lemma}
\proof Let $\mathcal{E}(x)$ denote the set of integers less than or
equal to $x$ with more than $(\log_2 x)^{\de/2}$ distinct prime
factors less than $y$ or more than $(\log_2 x)^{\de/2}$ distinct
prime factors greater than $z$. The size of this set is
\begin{align*}
|\mathcal{E}(x)| &\leq \frac{x}{\lfloor (\log_2 x)^{\de/2} \rfloor
!}\Big(\sum_{p\le y} \frac{1}{p} \Big)^{(\log_2 x)^{\de/2}}
+\frac{x}{\lfloor (\log_2 x)^{\de/2} \rfloor !} \Big( \sum_{z\le p\le x} \frac{1}{p} \Big) ^{(\log_2 x)^{\de/2}}\\
&\ll \frac{x}{(\log_2 x)^r},
\end{align*}
using that $\sum_{p\le x} \frac{1}{p}= \log_2 x + C +O(1/\log x)$
and Stirling's estimate $n!\sim \sqrt{2 \pi }n^{n+\half} e^{-n}$.
For $n\notin \mathcal{E}(x)$ we have $\om(n)-\om(n;y,z)\ll (\log_2
x)^{\de/2} $, and so it follows that
\begin{align*}
&\frac{\om(n)-\log_2 x}{\sqrt{\log_2 x}}-\frac{\om(n;y,z)-\sump
\frac{1}{p}}{\sqrt{\sump \frac{1}{p}}}\\
&=\frac{\om(n)-\log_2 x-\om(n;y,z)+\sump \frac{1}{p}}{\sqrt{\sump
\frac{1}{p}}}+O\Bigg(\frac{|\om(n)-\log_2 x|\log_3 x}{\log_2 x}\Bigg)\\
&\ll \frac{(\log_2 x)^{\de/2}}{\sqrt{\log_2 x}}+
\frac{|\om(n)-\log_2 x| \log_3 x}{\log_2 x}.
\end{align*}
The proof is complete by noting that except for $O\big(x/(\log_2
x)^r\big)$ integers less than $x$ we have $|\om(n)-\log_2 x|\ll
(\log_2 x)^{1/2+\de/4}$. This is because $\frac{\om(n)-\log_2
x}{\sqrt{\log_2 x}}$ has normal moments and in particular
$\sumn\Big(\frac{\om(n)-\log_2 x}{\sqrt{\log_2 x}}\Big)^{8\lfloor
r/\de \rfloor}\ll 1$, where the implied constant depends on $r/\de$.
\endproof

\section{Proof of Theorem \ref{fou}}

Define for a prime $p$
\begin{align*}
f_p(n)=
\begin{cases}
1-\frac{1}{p} &\text{ if } p|n \\
-\frac{1}{p} &\text{ if } p \nmid n \\
\end{cases},
\end{align*}
and if $m=\prod_i p_i ^ {\alpha_i}$, define $f_m(n)=\prod_i
f_{p_i}(n)^{\alpha_i}$. (Thus $f_1(n)=1$.) If we think of a prime
$p$ dividing $n$ with probability $1/p$ independently of other
primes, then we have $\mathbb{E}(f_m(n))=0$ for square-free $m$. So
we have written $\om(n)-\sum_{p\le x} \frac{1}{p}= \sum_{p\le x}
f_p(n)$ as a sum of independent random variables of mean $0$, which
already suggests by the Central Limit Theorem that
$\om(n)-\sum_{p\le x} \frac{1}{p}$ is normally distributed. This
simple idea is actually very powerful. It is borrowed from Granville
and Soundararajan \cite{gransound}, who use it to efficiently
compute very high moments of $\om(n)-\sum_{p\le x} \frac{1}{p}$ and
provide a new proof of the Erd{\H o}s-Kac theorem.

We have that the left hand side of Theorem \ref{fou} equals
\begin{align}
&\nonumber \sumn \prodi \exp\Big( \ii T_i \frac{\om(n+b_i;y,z)-\sump \frac{1}{p}}{\sqrt{\sump \frac{1}{p}}}\Big) \\
&\nonumber=\sumn \prodi \exp \Big(\sump \ii t_i f_p(n+b_i)\Big)  \\
&\nonumber=\sumn \prodp \exp \Big(\sumi \ii t_i f_p(n+b_i)\Big) \\
&\label{expan}=\sumn \prodp \Big(1+\Big(\sumi \ii t_i
f_p(n+b_i)\Big)+\frac{1}{2!}\Big(\sumi \ii t_i f_p(n+b_i)\Big)^2+
\ldots\Big).
\end{align}
Now upon expansion of the product, (\ref{expan}) equals
\begin{align}
\label{upon} \sumn \sum_{\substack{a_i \ge 1\\ p|a_i \Rightarrow
y<p<z}} K_{a_1,\ldots,a_r}  \prodi t_i^{\Omega(a_i)} f_{a_i}(n+b_i)
,
\end{align}
for some constants $K_{a_1,\ldots,a_r}$ of modulus bounded by 1,
where $\Omega(a)$ is the number of prime factors of $a$ counted with
multiplicity. Note that when the integers $a_i$ are pairwise coprime
we have that
\begin{align}
\label{k} K_{a_1,\ldots,a_r}=\prodi \prod_{p^{\al}|a_i}
\frac{\ii^{\al}}{\al!}.
\end{align}
We will evaluate (\ref{upon}) using the following results. The first
generalizes a result from \cite{gransound}.

\begin{lemma} \label{one}
Let $a_i$ be pairwise coprime integers for $1\le i \le r$. Denote
the square-free part of $a_i$ by $A_i$. We have
\begin{align*}
\sumn \prodi f_{a_i}(n+b_i)= \prodi & \prod_{p^{\alpha}
\parallel a_i} \Big( \frac{1}{p} \Big(1-\frac
{1}{p}\Big)^{\alpha} +\Big(\frac{-1}{p}\Big)^{\alpha} \Big(1-\frac
{1}{p}\Big)\Big)
\\
&+O\Big(\frac{1}{x}\prod_{i=1}^{r}\tau(A_i)^2\Big),
\end{align*}
where $\tau(A)$ denotes the number of divisors of $A$ and
$p^{\al}\parallel a$ means that $p^{\al}|a$ and $p^{\al+1}\nmid a$.
Note that the main term is zero unless each $a_i$ is square-full
(that is, $p|a_i$ implies $p^2|a_i$).
\end{lemma}

\proof For a fixed integer $a$ the value $f_a(n)$ depends only on
the common prime factors of $a$ and $n$, so $f_a(n)=f_a((A,n))$.
Thus we can group terms this way:
\begin{align}
\sumn \prodi f_{a_i}(n+b_i)
&\nonumber=\frac{1}{x}\sum_{d_i|A_i} \sum_{\substack{n\le x\\ (A_i, n+b_i)=d_i}} \prodi f_{a_i}(d_i) \\
&\nonumber=\frac{1}{x}\sum_{d_i|A_i} \sum_{\substack{n\le x\\ d_i|(n+b_i)}} \sum_{e_i|\frac{(A_i,n+b_i)}{d_i}} \prodi \mu(e_i) f_{a_i}(d_i) \\
&\label{one0}=\frac{1}{x}\sum_{d_i|A_i} \sum_{e_i|\frac{A_i}{d_i}}
\big(\sum_{\substack{n\le x\\ e_id_i|n+b_i}} 1\big) \prodi \mu(e_i)
f_{a_i}(d_i).
\end{align}
By the Chinese Remainder Theorem, since the integers $a_i$ are
pairwise coprime, $\sum_{\substack{n\le x\\ e_id_i|n+b_i}} 1=
\frac{x}{\prodi e_id_i} + O(1)$. Therefore the above sum is

\begin{align}
&\nonumber\prodi \sum_{d_i|A_i} \frac{f_{a_i}(d_i)}{d_i} \sum_{e_i|\frac{A_i}{d_i}} \frac{\mu(e_i)}{e_i} +O\Big(\frac{1}{x}\prod_{i=1}^{r}\tau(A_i)^2\Big) \\
&\nonumber=\prodi \sum_{d_i|A_i} \frac{f_{a_i}(d_i)}{d_i} \frac{\phi(\frac{A_i}{d_i})}{\frac{A_i}{d_i}} +O\Big(\frac{1}{x}\prod_{i=1}^{r}\tau(A_i)^2\Big)  \\
&=\prodi \sum_{d_i|A_i} f_{a_i}(d_i) \frac{1}{A_i}
\phi\Big(\frac{A_i}{d_i}\Big)
+O\Big(\frac{1}{x}\prod_{i=1}^{r}\tau(A_i)^2\Big).
\end{align}
Now it is easily verified (by multiplicativity in $a_i$) that the
main term of the last line above equals
\begin{align*}
\prodi \prod_{p^{\alpha}
\parallel a_i} \Big( \frac{1}{p} \Big(1-\frac
{1}{p}\Big)^{\alpha} +\Big(\frac{-1}{p}\Big)^{\alpha} \Big(1-\frac
{1}{p}\Big)\Big).
\end{align*}
\endproof
\noindent In the case that $a_1,\ldots,a_r$ are not pairwise coprime
we will need the following.

\begin{lemma} \label{two} Let $r\ge 2$ and $0\le b_1<...<b_r\le
\lambda (\log_2 x)^{1/2-\de}$. Suppose that for some prime $y<q<z$,
we have that $q|a_1$ and $q|a_2$. Let $q^{\al_i}
\parallel a_i$ and let $a_i'=a_i q^{-\al_i}$. Let $A_i$ denote the squarefree part of $a_i$ and let $A_i'$ denote the square free part of $a_i'$. We have
\begin{align*}
\sumn \prodi f_{a_i}(n+b_i) = O\Big(\frac{1}{q^2}\sumn \prodi
f_{a_i'}(n+b_i)\Big)
+O\Big(\frac{1}{qx}\prod_{i=1}^{r}\tau(A_i)^2\Big).
\end{align*}
\end{lemma}
\proof From (\ref{one0}) we have
\begin{align}
\label{saw} \sumn \prodi f_{a_i}(n+b_i) = \frac{1}{x}\sum_{d_i|A_i}
\sum_{e_i|\frac{A_i}{d_i}} \big(\sum_{\substack{n\le x\\
e_id_i|n+b_i}} 1\big) \prodi \mu(e_i) f_{a_i}(d_i).
\end{align}
The sum $\sum_{\substack{n\le x\\ e_id_i|n+b_i}} 1$ is zero for
large enough $x$ if $q$ divides more than one integer $e_i d_i$.
This is because if $q|n+b_i$ and $q|n+b_j$ then $q|b_i-b_j$ and
hence $i=j$ since $q> y$ and the integers $b_i$ are distinct and
bounded by $\la (\log_2 x)^{1/2-\de}$. Thus we can suppose that $q$
divides at most one integer $e_i d_i$. First consider the terms of
(\ref{saw}) with $q \nmid e_id_i$ for all $i$. These terms
contribute
\begin{align*}
&\frac{1}{x}\sum_{d_i|A_i'} \sum_{e_i|A_i'/d_i} \big(\sum_{\substack{n\le x\\ e_id_i|n+b_i}} 1\big) \prodi \mu(e_i) f_{a_i'}(d_i)f_{q^{\al_i}}(d_i)\\
&=\Big(\frac{-1}{q}\Big)^{\alpha_1+\ldots+\alpha_r} \sumn \prodi
f_{a_i'}(n+b_i),
\end{align*}
using that $f_{q^{\al_i}}(d_i)=(-1/q)^{\al_i}$ and the identity of
(\ref{saw}). Since $\al_1+\al_2 \ge 2$ we get the desired factor of
$1/q^2$. Now say $q|e_1d_1$ and $q\nmid e_id_i$ for $i\neq 1$. The
contribution of this case is,
\begin{align}
&\label{two1}\frac{1}{x}\sum_{d_i|A_i'} \sum_{e_i|A_i'/d_i} \big(\sum_{\substack{n\le x\\ qe_1d_1|n+b_1\\ e_id_i|n+b_i,i\neq 1}} 1\big) \mu(e_1) f_{a_1'}(qd_1)f_{q^{\al_1}}(qd_1) \prod_{i=2}^{r} \mu(e_i) f_{a_i'}(d_i)f_{q^{\al_i}}(d_i) \\
&\nonumber+\frac{1}{x}\sum_{d_i|A_i'} \sum_{e_i|A_i'/d_i'}
\big(\sum_{\substack{n\le x\\ qe_1d_1|n+b_1\\ e_id_i|n+b_i,i\neq 1}}
1\big) \mu(qe_1) f_{a_1'}(d_1)f_{q^{\al_1}}(d_1) \prod_{i=2}^{r}
\mu(e_i) f_{a_i'}(d_i)f_{q^{\al_i}}(d_i) ,
\end{align}
where the first line corresponds to $q|d_1$ and the second to
$q|e_1$. We have by the Chinese Remainder Theorem,
\begin{align*}
\sum_{\substack{n\le x\\
qe_1d_1|n+b_1\\e_id_i|n+b_i,i\neq 1}} 1
=\frac{1}{q}\sum_{\substack{n\le x\\e_id_i|n+b_i}} 1 +O(1).
\end{align*}
Thus the contribution of the sums of (\ref{two1}) is
\begin{align}
&\label{two2}\Big(\frac{1}{q}\Big(1-\frac{1}{q}\Big)^{\al_1} \Big(\frac{-1}{q}\Big)^{\alpha_2+\ldots+\alpha_r} - \frac{1}{q}\Big(\frac{-1}{q}\Big)^{\alpha_1+\ldots+\alpha_r}\Big)  \sumn \prodi f_{a_i'}(n+b_i) \\
&\nonumber\indent +O\Big(\frac{1}{x}|f_{q^{\al_2}}(d_2)| \prodi
\sum_{e_id_i|A_i'} 1\Big).
\end{align}
Again we have a factor of $1/q^2$ in the first line above since
$\al_2\ge 1$. The second line of (\ref{two2}) is
$O(\frac{1}{qx}\prod_{i=1}^{r}\tau(A_i)^2)$ since $\al_2\ge 1$ and
$q\nmid d_2$. This completes the proof as terms with $q|e_jd_j$ and
$q\nmid e_id_i$ for $i\neq j$ are dealt with similarly.
\endproof
\noindent We will also need the following observations.

\begin{lemma} \label{three} We have
\begin{align*}
\frac{1}{x}\sum_{n\le x} \sum_{\substack{a_i \ge 1 \\ p|a_i \Rightarrow y<p<z\\
\omega(a_1)\ge (\log_2 x)^{2r}}} |K_{a_1,\ldots,a_r}|
\prod_{i=1}^{r} |f_{a_i}(n+b_i)| |t_i|^{\Omega(a_i)}\ll
\frac{1}{\log x}.
\end{align*}
\end{lemma}
\proof We first bound the contribution of terms with
$\omega(a_i)=w_i$ for some positive integers $w_i$. Recall that
$|K_{a_1,\ldots,a_r}|\le 1$ and note that $|f_{p^{\alpha}}(n)|\le
|f_p(n)|$. Thus
\begin{align}
&\label{three1} \frac{1}{x}\sum_{n\le x} \sum_{\substack{a_i \ge 1 \\ p|a_i \Rightarrow y<p<z\\
\omega(a_i)=w_i}} |K_{a_1,\ldots,a_r}| \prod_{i=1}^{r}
|f_{a_i}(n+b_i)|
|t_i|^{\Omega(a_i)}\\
&\nonumber \ll \frac{1}{x}\sum_{n\le x} \sum_{\substack{a_i \ge 1 \\ p|a_i \Rightarrow y<p<z\\
\omega(a_i)=w_i}}  \prod_{i=1}^{r} |f_{A_i}(n+b_i)|
|t_i|^{\Omega(a_i)}
\end{align}
For a fixed square-free integer $A_i$ with $\omega(A_i)=w_i$ we have
\begin{align*}
\sum_{\substack{a_i\ge 1\\ \text{$A_i$=square-free part of $a_i$}}}
|t_i|^{\Omega(a_i)} \le 1,
\end{align*}
since $|t_i|\le \frac{1}{1000}$. Thus (\ref{three1}) is bounded by
\begin{align}
\label{three2}\ll \frac{1}{x}\sum_{n\le x} \prod_{i=1}^{r}
\frac{1}{w_i!}\Big(\sum_{y<p<z} |f_p(n+b_i)|\Big)^{w_i}.
\end{align}
Since $|f_p(n)|\le 1$ if $p|n$ and $|f_p(n)|\le \frac{1}{p}$ if
$p\nmid n$, we have that (\ref{three2}) is bounded by
\begin{align}
&\label{three3}\ll \frac{1}{x}\sum_{n\le x} \prod_{i=1}^{r}
\frac{1}{w_i!}\Big(\omega(n+b_i;y,z)+\log_2 x\Big)^{w_i}\\
&\nonumber \ll \Big(\prod_{i=1}^{r} \frac{1}{w_i!}\Big)
\frac{1}{x}\sum_{n\le x} \sum_{i=1}^{r}
2^{rw_i}\Big(\omega(n+b_i;y,z)^{rw_i}+ (\log_2 x)^{rw_i}\Big).
\end{align}
We have
\begin{align}
\label{three4}\frac{1}{x}\sum_{n\le x} \omega(n+b_i;y,z)^{rw_i}\ll
\frac{1}{x} \sum_{y<p_{1},\ldots,p_{rw_i}<z}\sum_{\substack{n\le x\\
[p_{1},\ldots, p_{rw_i}]|n+b_i}} 1,
\end{align}
where $[p_{1},\ldots, p_{rw_i}]$ denotes the lowest common multiple
of $p_{1},\ldots, p_{rw_i}$. Now (\ref{three4}) is bounded by
\begin{align}
\ll
\sum_{y<p_{1},\ldots,p_{rw_i}<z}\frac{1}{[p_{1},\ldots,p_{rw_i}]}\ll
2^{rw_i}(\log_2 x)^{rw_i}.
\end{align}
Thus we have that (\ref{three3}) is bounded by
\begin{align}
\ll\Big( \prod_{i=1}^{r} \frac{1}{w_i!}\Big) \sum_{i=1}^{r} 4^{rw_i}
(\log_2 x)^{rw_i}.
\end{align}
Summing over integers $w_i\ge 1$ for $i\ge 2$ this is bounded by
\begin{align}
\label{three5}\ll \frac{(4\log_2 x)^{rw_1}\exp( (4\log_2 x)^r
)}{w_1!} .
\end{align}
Finally the sum of (\ref{three5}) over integers $w_1\ge (\log_2
x)^{2r}$ is $\ll 1/\log x$.
\endproof

\begin{lemma} \label{four}
We have
\begin{align*}
 \sum_{\substack{a_i\ge 1\\ p|a_i\Rightarrow y<p<z\\ \omega(a_1)\ge
(\log_2 x)^{2r}}} \prod_{i=1}^{r} |t_i|^{\Omega(a_i)}
\prod_{p^{\alpha}\parallel a_i} \Big|
\frac{1}{p}\Big(1-\frac{1}{p}\Big)^{\alpha} +
\Big(1-\frac{1}{p}\Big)\Big(\frac{-1}{p}\Big)^{\alpha}\Big|\ll
\frac{1}{\log x}.
\end{align*}
\end{lemma}
\proof We first bound the contribution of terms with
$\omega(a_i)=w_i$ for some positive integers $w_i$. We have
\begin{align}
&\label{four1} \sum_{\substack{a_i\ge 1\\ p|a_i\Rightarrow y<p<z\\
\omega(a_i)=w_i }} \prod_{i=1}^{r} |t_i|^{\Omega(a_i)}
\prod_{p^{\alpha}\parallel a_i} \Big|
\frac{1}{p}\Big(1-\frac{1}{p}\Big)^{\alpha} +
\Big(1-\frac{1}{p}\Big)\Big(\frac{-1}{p}\Big)^{\alpha}\Big|\\
&\nonumber \ll \prod_{i=1}^{r} \frac{1}{w_i!} \Big(\sum_{y<p<z}
\frac{1}{p} \Big)^{w_i}\ll \prod_{i=1}^{r}  \frac{(\log_2
x)^{w_i}}{w_i!}.
\end{align}
Summing over integers $w_i\ge 1$ for $i\ge 2$ this is bounded by
\begin{align}
\label{four2} \ll \frac{(\log_2 x)^{w_1}\exp(r\log_2 x)}{w_1!} .
\end{align}
The sum of (\ref{three5}) over integers $w_1\ge (\log_2 x)^{2r}$ is
$\ll 1/\log x$.
\endproof

\subsection*{Back to the proof}

By Lemma \ref{three} we see that (\ref{upon}) equals, up to an error
of $O(1/\log x)$,
\begin{align}
\label{upon2} \frac{1}{x}\sum_{n\le x} \sum_{\substack{a_i \ge
1\\p|a_i \Rightarrow y<p<z\\ \omega(a_i)\le (\log_2 x)^{2r}}}
K_{a_1,\ldots,a_r} \prod_{i=1}^{r} f_{a_i}(n+b_i) t_i^{\Omega(a_i)}.
\end{align}
Let us first treat the terms of (\ref{upon2}) with $a_1,\ldots,a_r$
not pairwise coprime. Applying Lemma \ref{two} repeatedly, these
terms contribute an amount bounded by
\begin{align}
&\label{notcoprimeerror} \ll\sum_{y<q<z} \frac{1}{q^2} \sum_{\substack{a_i \text{ pairwise coprime}\\
p|a_i \Rightarrow y<p<z\\ \omega(a_i)\le (\log_2 x)^{2r} }}
\Big|\frac{1}{x}\sum_{n\le x} \prod_{i=1}^{r}
f_{a_i}(n+b_i) t_i^{\Omega(a_i)}\Big|\\
&\nonumber+\sum_{y<q<z} \frac{1}{qx} \sum_{\substack{a_i\ge 1\\
p|a_i \Rightarrow y<p<z\\ \omega(a_i)\le (\log_2 x)^{2r}}}
\prod_{i=1}^{r} |t_i|^{\Omega(a_i)} \tau(A_i)^2.
\end{align}
For the second line of (\ref{notcoprimeerror}), we have
\begin{align}
&\label{already} \frac{1}{x}\prod_{i=1}^{r} \sum_{\substack{a_i\ge 1\\
p|a_i \Rightarrow y<p<z\\ \omega(a_i)\le (\log_2 x)^{2r}}}
|t_i|^{\Omega(a_i)}
\tau(A_i)^2\ll \frac{1}{x}\prod_{i=1}^{r} \sum_{\substack{a_i \text{square-free} \\
p|a_i \Rightarrow y<p<z\\ \omega(a_i)\le (\log_2 x)^{2r}}} \tau(A_i)^2\\
&\nonumber \ll \frac{1}{x}\prod_{i=1}^{r} z^{(\log_2 x)^{2r}}
4^{(\log_2 x)^{2r}} \ll x^{-1/2}.
\end{align}
Thus the second line of (\ref{notcoprimeerror}) falls into the error
term of Theorem \ref{fou}. To bound the first line of
(\ref{notcoprimeerror}) we use Lemma \ref{one} to get that
\begin{align}
&\sum_{y<q<z} \frac{1}{q^2} \sum_{\substack{a_i \text{ pairwise coprime}\\
p|a_i \Rightarrow y<p<z\\ \omega(a_i)\le (\log_2 x)^{2r}}}
\Big|\frac{1}{x}\sum_{n\le x} \prod_{i=1}^{r}
f_{a_i}(n+b_i) t_i^{\Omega(a_i)}\Big|\\
&\nonumber \ll\frac{1}{y} \Big(\sum_{\substack{a_i\ge 1\\
p|a_i \Rightarrow y<p<z\\ \omega(a_i)\le (\log_2 x)^{2r}}}
\prod_{i=1}^{r} |t_i|^{\Omega(a_i)}\prod_{p^{\alpha}\parallel a_i}
\frac{1}{p} + \sum_{\substack{a_i\ge 1\\
p|a_i \Rightarrow y<p<z\\ \omega(a_i)\le (\log_2 x)^{2r}}}
\frac{1}{x} \prod_{i=1}^{r}
|t_i|^{\Omega(a_i)} \tau(A_i)^2 \Big)\\
&\nonumber \ll \frac{1}{y} \prod_{i=1}^{r} \prod_{y<p<z}
\Big(1+\frac{1}{p}\Big)+x^{-1/2},
\end{align}
where we used the bound of (\ref{already}). Now this is less than $
\frac{ (\log x)^r }{y}+x^{-1/2}\ll 1/\log x$.

Thus only the terms of (\ref{upon2}) with $a_1,\ldots,a_r$ coprime
will give a main contribution. Using Lemma \ref{one} and (\ref{k})
we get
\begin{align}
&\label{upon3} \frac{1}{x}\sum_{n\le x} \sum_{\substack{a_i \text{ pairwise coprime}\\  p|a_i \Rightarrow y<p<z\\
\omega(a_i)\le (\log_2 x)^{2r}}}  K_{a_1,\ldots,a_r} \prod_{i=1}^{r}
f_{a_i}(n+b_i)
t_i^{\Omega(a_i)}\\
=&\nonumber \prod_{i=1}^{r} \sum_{\substack{a_i \text{ pairwise coprime}\\ p|a_i\Rightarrow y<p<z\\
\omega(a_i)\le (\log_2 x)^{2r}}} \prod_{p^{\alpha}\parallel
a_i}\frac{{\bf i}^{\alpha}}{\alpha !} t_i^{\alpha} \Big(
\frac{1}{p}\Big(1-\frac{1}{p}\Big)^{\alpha} +
\Big(1-\frac{1}{p}\Big)\Big(\frac{-1}{p}\Big)^{\alpha}\Big)\\
&\nonumber\indent+O\Big( \sum_{\substack{a_i\ge 1\\ p|a_i
\Rightarrow y<p<z\\ \omega(a_i)\le (\log_2 x)^{2r}}}
\frac{1}{x}\prod_{i=1}^{r} |t_i|^{\Omega(a_i)} \tau(A_i)^2\Big).
\end{align}
We've already seen in (\ref{already}) that the error term above is
negligible. The main term of (\ref{upon3}) is zero if $a_i$ is not
square-full for all $i$. Thus we may further impose the condition
$p|a_i \Rightarrow p^2|a_i$. We may also extend by Lemma \ref{four}
the sum in the main term of (\ref{upon3}) to {\it all} pairwise
coprime and square-full integers $a_i$ whose prime factors lie
between $y$ and $z$, up to an error of $O(1/\log x)$. Thus
(\ref{upon3}) equals up to this error,
\begin{align}
&\label{final1} \prod_{i=1}^{r} \prod_{y<p<z} \Big( 1+ \frac{1}{p}
\sum_{\alpha\ge 2} \frac{{\bf i}^{\alpha}}{\alpha !} t_i^{\alpha}
\Big(1-\frac{1}{p}\Big)^{\alpha} + \Big(1-\frac{1}{p}\Big)
\sum_{\alpha \ge 2} \frac{{\bf i}^{\alpha}}{\alpha !}
t_i^{\alpha} \Big(-\frac{1}{p}\Big)^{\alpha} \Big)\\
&=\nonumber \prod_{i=1}^{r} \prod_{y<p<z} \Big(1+ \frac{1}{p}\Big(
e^{{\bf i} t_i}-1-{\bf i} t_i
\Big)+O\Big(\frac{1}{p^{2}}\Big) \Big)\\
&\nonumber =\prod_{i=1}^{r} \exp\Big( \sum_{y<p<z} \log \Big(1+
\frac{1}{p}\Big( e^{{\bf i} t_i}-1-{\bf i} t_i
\Big)+O\Big(\frac{1}{p^{2}}\Big)\Big)\Big).
\end{align}
Now since
\begin{align*}
\log \Big(1+ \frac{1}{p}\Big( e^{{\bf i} t_i}-1-{\bf i} t_i
\Big)+O\Big(\frac{1}{p^{2}}\Big)\Big)=\frac{1}{p}\Big( e^{{\bf i}
t_i}-1-{\bf i} t_i \Big)+O\Big(\frac{1}{p^{2}}\Big),
\end{align*}
we have that (\ref{final1}) equals
\begin{align}
&\prod_{i=1}^{r} \exp\Big( \Big( e^{{\bf i} t_i}-1-{\bf i} t_i \Big)
\sum_{y<p<z}
\frac{1}{p}+O\Big(\frac{1}{y}\Big)\Big)\\
&\nonumber=\prod_{i=1}^{r} \exp\Big( \Big( e^{{\bf i} t_i}-1-{\bf i}
t_i \Big) \sum_{y<p<z} \frac{1}{p}\Big)+O(1/\log x).
\end{align}

{\bf Acknowledgments.} I am grateful to Prof. K. Soundararajan for
posing the problem studied in this paper and for many helpful
discussions pertaining to it. I am also thankful to Prof. A.
Granville for numerous useful comments. Part of this work was done
with support through a grant from the NSF (DMS 0500711) and while I
visited the Centre de Recherches Math$\acute{\text{e}}$matiques,
Montr$\acute{\text{e}}$al.

\nocite{sathe} \nocite{selberg} \nocite{turan} \nocite{feller}

\bibliographystyle{amsplain}

\bibliography{poisson}

\end{document}